\documentstyle[12pt]{article}
\newcommand{\be}{\begin{equation}}
\newcommand{\ee}{\end{equation}}

\newcommand{\al}{\alpha}

\newcommand{\bet}{\beta}

\newcommand{\M}{{\cal M}}

\newcommand{\eps}{\varepsilon}

\newcommand{\dz}{\wedge}
\newcommand{\C}{{\bf C }}

\newcommand{\ba}{\begin{array}}
\newcommand{\ea}{\end{array}}
\newcommand{\beq}{\begin{eqnarray}}
\newcommand{\eeq}{\end{eqnarray}}
\newtheorem{lm}{Lemma}
\newtheorem{th}{Theorem}
\newtheorem{pr}{Proposition}
\newtheorem{co}{Corollary}
\newtheorem{rem}{Remark}
\newtheorem{deff}{Definition}
\newcommand{\bd}{\begin{deff}}
\newcommand{\ed}{\end{deff}}
\newcommand{\bl}{\begin{lm}}
\newcommand{\el}{\end{lm}}
\newcommand{\bp}{\begin{pr}}
\newcommand{\ep}{\end{pr}}
\newcommand{\bt}{\begin{th}}
\newcommand{\et}{\end{th}}
\newcommand{\bc}{\begin{co}}
\newcommand{\ec}{\end{co}}
\newcommand{\brm}{\begin{rem}}
\newcommand{\erm}{\end{rem}}

\newcommand{\der}{{\rm d}}
\begin{document}

\thispagestyle{empty}

\title {A FOUR DIMENSIONAL EXAMPLE OF RICCI FLAT METRIC ADMITTING   
ALMOST-K\"AHLER NON-K\"AHLER STRUCTURE
\footnote{Research supported in part by: Komitet Bada\'n Naukowych 
(Grant nr 2 P302 112 7), Consorzio per lo Sviluppo Internazionale 
dell'Universita degli Studi di Trieste and 
Erwin Schr\"{o}dinger 
International Institute for Mathematical Physics.}\\ 
\vskip 1.truecm
{\small {\sc Pawe\l  ~Nurowski}}\\
{\small {\it Dipartimento di Scienze Matematiche }}\\
\vskip -0.3truecm
{\small {\it Universita degli Studi di Trieste, Trieste, Italy}}
\footnote{Permanent address: {\it Katedra Metod Matematycznych Fizyki, 
Wydzia\l ~Fizyki, Uniwersytet Warszawski, ul. Ho\.za 74, Warszawa, 
Poland, e-mail: nurowski@fuw.edu.pl}}\\
{\small {\sc Maciej Przanowski}}\\
{\small {\it Instytut Fizyki, Politechnika \L\'{o}dzka }}\\
\vskip -0.3truecm
{\small {\it W\'{o}lcza\'nska 219, 93-005 \L\'{o}d\'z, Poland}}}
\author{\mbox{}}
\maketitle
\begin{abstract}
\noindent
We construct an example of Ricci-flat almost-K\"ahler non-K\"ahler 
structure in four dimensions.
\end{abstract}
\newpage
\noindent
1. Let $\cal M$ be a 4-manifold equipped with a metric $g$ of signature 
(++++). The pair ($\cal M$, $g$) is called a Riemmanian 4-manifold.\\

\noindent
An almost hermitian structre on ($\cal M$, $g$) is a tensor field 
$J: T\M\to T\M$ such that $J^2=-id$ and $g(JX,JY)=g(X,Y)$. An almost 
hermitian structure ($\cal M$, $g$, $J$) is called hermitian if $J$ is 
integrable. Due to the Newlander-Nirenberg theorem this is equivalent 
to the vanishing of the Nijenhuis tensor 
$N_J(X,Y)=[JX,JY]-[X,Y]-J[JX,Y]-J[X,JY]$ for $J$.\\ 

\noindent
Given an almost hermitian structure ($\cal M$, $g$, $J$) one defines 
the fundamental 2-form $\omega$ by $\omega(X,Y)=g(X,JY)$. An almost hermitian 
structure ($\cal M$, $g$, $J$) is called almost-K\"ahler if its fundamental 
2-form is closed. If, in addition, $J$ is integrable then such structure 
is called K\"ahler.\\

\noindent
This paper is motivated by the following conjecture \cite{bi:Goldberg}.\\

\noindent
{\bf Goldberg's Conjecture}\\
{\it The almost K\"{a}hler structure of a compact Einstein manifold is 
necessarilly K\"{a}hler.}\\

\noindent
The conjecture was proven in the case of non-negative scalar curvature of the 
Einstein manifold by K. Sekigawa in \cite{bi:Sekigawa}.\\

\noindent
In this paper we show that the assumption about compactness of the Einstein 
manifold is essential for the Goldberg conjecture. In particular, we give an 
explicit example of a Ricci-flat almost-K\"ahler 
non-K\"ahler structure on a noncompact 4-manifold. 
This result is given by Theorem 1 of paragraph 4.\\

\noindent
2. Let $\cal U$ be an open subset of ${\bf R}^4$. Let 
$\theta^i=(M,\bar{M}, N,\bar{N})$ be four complex-valued 1-forms on 
$\cal U$  
such that $M\dz\bar{M}\dz N\dz\bar{N}\neq 0$. Using $\theta^i$ we define a 
metric $g$ on $\cal U$ by 
$$
g=2(M\bar{M}+ N\bar{N}):=M\otimes\bar{M}+\bar{M}\otimes M+ N\otimes\bar{N}+ 
\bar{N}\otimes N.
$$
Clearly ($\cal U$, $g$) is a Riemannian 4-manifold.\\
The Weyl tensor $W$ of the metric $g$ splits onto self-dual $(W^+)$ and 
anti-self-dual $(W^-)$ parts. $({\cal U},g)$ is said to be (anti-)self-dual 
iff 
($W^+\equiv 0$) $W^-\equiv 0$. If ($W^+\neq 0$) $W^-\neq 0$ then in every 
point of 
$\cal U$ it defines at most two spinor directions ($[\al^+,\bet^+]$) 
$[\al^-,\bet^-]$; see e.g. \cite{bi:phd,bi:przan1}. 
$(W^+)$ $W^-$ is said to be of type 
$D$ if $(\al^+)$ $\al^-$ coincides with ($\bet^+$) $\bet^-$.\\ 

\noindent
Let $e_i= (m,\bar{m}, n,\bar{n})$ be a basis dual to 
$\theta^i=(M,\bar{M}, N,\bar{N})$. For any $\xi\in\C$ it is convenient to 
consider 1-forms  
$$
M_\xi=\frac{M-\bar{\xi}\bar{N}}{\sqrt{1+\xi\bar{\xi}}}\quad\quad
N_\xi=\frac{N+\bar{\xi}\bar{M}}{\sqrt{1+\xi\bar{\xi}}}
$$
and vector fields
$$
m_\xi=\frac{m-\xi\bar{n}}{\sqrt{1+\xi\bar{\xi}}}\quad\quad
n_\xi=\frac{n+\xi\bar{m}}{\sqrt{1+\xi\bar{\xi}}}.
$$
The following Lemma is well known (see for example \cite{bi:phd,bi:przan1}).
\bl~\\
\vskip -1.truecm
\begin{itemize}
\item[{\bf i)}] 
For any value of the complex parameter $\xi\in\C\cup\{\infty\}$
the expressions 
$$
J^+_\xi=i(\overline{M_\xi}\otimes\overline{m_\xi}-M_\xi\otimes m_\xi+
\overline{N_\xi}\otimes\overline{n_\xi}-N_\xi\otimes n_\xi)
$$
$$
J^-_\xi=i(M_\xi\otimes m_\xi-\overline{M_\xi}\otimes\overline{m_\xi}+
\overline{N_\xi}\otimes\overline{n_\xi}-N_\xi\otimes n_\xi)
$$
define almost hermitian structures on $({\cal U}, g)$. 
\item[{\bf ii)}] 
The fundamental 2-forms corresponding to $J^+_\xi$ and $J^-_\xi$ are 
respectively given by 
$$
\omega^+_\xi=i(M_\xi\dz\overline{M_\xi}+N_\xi\dz
\overline{N_\xi})
$$
$$
\omega^-_\xi=i(\overline{M_\xi}\dz M_\xi+N_\xi\dz\overline{N_\xi}).
$$
\item[{\bf iii)}] 
Any almost hermitian structure on $({\cal U}, g)$ is given either by one of 
$J^+_\xi$ or by one of $J^-_\xi$. Structures $J^+_\xi$ are different from 
$J^-_\xi$; also, different $\xi$s correspond to different structures.
\item[{\bf iv)}]
If the metric $g$ is not 
self-dual then among $J^+_\xi$s only at most four structures, 
corresponding to  
specific four values of the parameter $\xi$, may be integrable. 
Analogously, if the metric $g$ is not anti-self-dual then 
only at most four $J^-_\xi$s may be integrable.
\end{itemize}
\el

\noindent
3. Let $(x^1,x^2,x^3,x^4)$ be Euclidean coordinates on $\cal U$. Define 
\be
z_1=x^1+ix^2 \quad\quad\quad z_2=x^3+ix^4. 
\label{eq:coo}
\ee
Let $\partial_k=\frac{\partial}{\partial z_k}$ and 
$\partial_{\bar{k}}=\frac{\partial}{\partial\overline{z_k}}$, $k=1,2$.\\  

\noindent
Consider two 1-forms $M$ and $N$ on $\cal U$ defined by
\be
M=f(\der z_1+h\der z_2)\quad\quad\quad N=\frac{1}{f}\der z_2,
\label{eq:mn}
\ee
where $f\neq 0$ (real) and $h$ (complex) are functions on $\cal U$.\\

\noindent
Since $M\dz\bar{M}\dz N\dz\bar{N}=\der z_1\dz\der\overline{z_1}\dz\der 
z_2\dz\der\overline{z_2}\neq 0$ then the metric $g=2(M\bar{M}+N\bar{N})$ 
equipes $\cal U$ with the Riemannian structure. Consider  
almost hermitian structures $J^+_\xi$ for such ($\cal U$, $g$). It is 
interesting to note that if $\xi={\rm e}^{i\phi}=$const then the 
corresponding fundamental 2-form $\omega^+_{{\rm e}^{i\phi}}$ reads 
$$
\omega^+_{{\rm e}^{i\phi}}=i({\rm e}^{i\phi}\der z_2\dz\der z_1-
{\rm e}^{-i\phi}\der\overline{z_2}\dz\der\overline{z_1}) 
$$
and is closed. Thus, for any ${\rm e}^{i\phi}\in{\bf S}^1$ we constructed 
an almost-K\"ahler structure ($\cal U$, $g$, $J^+_{{\rm e}^{i\phi}}$). 
If the functions $f$ and $h$ are general enough, then the 
metric $g$ has no chance to be self-dual. Moreover, since in such case  
there is a finite number of hermitian structures among $J^+_\xi$, then most 
of our structures must be non-K\"ahler. Summing up we have the folowing 
Lemma.
\bl
Let $(z_1,\overline{z_1}, z_2,\overline{z_2})$ be coordinates on $\cal U$ as 
in (\ref{eq:coo}). Then for each value of the real constant 
$\phi\in [0,2\pi[$ the metric  
\be
g=2f^2(\der z_1+h\der z_2)(\der\overline{z_1}+\bar{h}\der\overline{z_2})+
2\frac{1}{f^2}\der z_2\der\overline{z_2}
\label{eq:mmm}
\ee
and the almost complex structure
\be
J^+_{{\rm e}^{i\phi}}=2{\rm Re}\{
i{\rm e}^{i\phi}[f^2(\der z_1+h\der z_2)\otimes 
(\partial_{\bar{2}}-\bar{h}\partial_{\bar{1}})-\frac{1}{f^2}\der z_2\otimes
\partial_{\bar{1}}]\}
\label{eq:mm}
\ee
define an almost-K\"ahler structure on $\cal U$.\\
If the functions $f$ and $h$  
are general enough to prevent the metric of beeing self-dual then 
these structures are non-K\"ahler for almost all values of $\phi$.
\el

\noindent
4. We look for not-self-dual Ricci-flat metrics among the metrics 
of Lemma 2. For this pourpose it is convenient to restrict to the  
metrics (\ref{eq:mmm}) whose anti-self-dual part of the 
Weyl tensor is strictly of 
type D. Such a restriction guarantees that all structures 
(\ref{eq:mm}) are non-K\"ahler \cite{bi:phd,bi:przan1}.\\

\noindent
We recall a useful Lemma \cite{bi:przan3}.\\
\bl
Let $g$ be a Ricci-flat Riemannian metric in four dimensions. Assume that 
the anti-self-dual part of the 
Weyl tensor for $g$ is strictly of type D. Then, 
locally there always exist 
complex coordinates $(z_1, z_2)$ and a real function 
$K=K(v,z_2,\overline{z_2})$, $v=z_1+\overline{z_1}$ such that the metric can   
can be written as
\be
g=\frac{\eps K_{vv}}{(K_v)^{3/2}}
(\der z_1+\frac{K_{v{\small 2}}}{K_{vv}}\der z_2)
(\der \overline{z_1} +\frac{K_{v\small\bar{2}}}{K_{vv}}\der 
\overline{z_2})+4{\rm e}^{-K}\frac{(K_v)^{1/2}}{\eps K_{vv}}\der z_2\der
\overline{z_2},
\label{eq:przm}
\ee
where $K_{v\small\bar{2}}=\frac{\partial^2 K}{\partial v\partial
\overline{z_2}}$, etc. The function $K$ satisfies 
\be
K_{vv}K_{\small 2\bar{2}}-K_{v\small\bar{2}}K_{v\small 2}-2
{\rm e}^{-K}(K_{vv}+2(K_v)^2)=0,
\label{eq:k}
\ee
\be
K_v>0,\quad\quad\quad\eps K_{vv}>0
\label{eq:kk}
\ee
where $\eps$ is either plus or minus one.\\
Also, every function $K=K(v,z_2,\overline{z_2})$ satisfying 
(\ref{eq:k})-(\ref{eq:kk}) defines, via (\ref{eq:przm}), a Ricci-flat metric. 
This metric has the anti-self-dual part of the Weyl tensor of strictly type D.
\el

\noindent
We ask when the metric (\ref{eq:mmm}) can be written in the form 
(\ref{eq:przm}). Identifying coordinates $(z_1,z_2)$ in both metrics 
we see that it is possible if 
$$2f^2=\frac{\eps K_{vv}}{(K_v)^{3/2}}\quad\quad\quad {\rm and} 
\quad\quad\quad  
\frac{2}{f^2}=4{\rm e}^{-K}\frac{(K_v)^{1/2}}{\eps K_{vv}}.$$
These two equations are compatible only if $K_v{\rm e}^K=1$. 
It is a matter of straightforward integration that, modulo the coordinate 
transformations, the general solution of 
this equation which simultaneously satisfies the equation (\ref{eq:k}) 
is\footnote{This solution was already known to S\l awomir Bia\l ecki in 1984 
\cite{bi:bia}.} 
$K=\log (v-2z_2\overline{z_2})$. Using such $K$ we easily 
find that in the region 
$$
{\cal U}'=\{{\cal U}\ni (z_1,z_2)\quad{\rm s.t.}\quad 
v-2z_2\overline{z_2} >0\}
$$
the metric (\ref{eq:mmm}) with 
$$
f=\frac{1}{\sqrt{2}(v-2z_2\overline{z_2})^{1/4}},
\quad\quad\quad h=-2\overline{z_2},
$$
is Ricci-flat and strictly of type D on the anti-self-dual side of its Weyl 
tensor. The explicit expression for such $g$ reads
\be
g=\frac{1}{(v-2z_2\overline{z_2})^{1/2}}
(\der z_1-2\overline{z_2}\der z_2)(\der\overline{z_1}-
2z_2\der\overline{z_2})+
4(v-2z_2\overline{z_2})^{1/2}\der z_2\der\overline{z_2},
\label{eq:metr}
\ee

\noindent
To have a better insight into this 
metric we choose new coordinates 
$$
x=(v-2z_2\overline{z_2})^{1/2},\quad\quad\quad y=z_2+\overline{z_2},
\quad\quad\quad z=i(\overline{z_2}-z_2),\quad\quad\quad 
q=\frac{z_1-\overline{z_1}}{2i}
$$ 
on ${\cal U}'$. 
These coordinates are real. The metric (\ref{eq:metr}) in 
these coordinates reads
$$
g=x(\der x^2+\der y^2+\der z^2)+\frac{1}{x}
(\frac{1}{2}z\der y-\frac{1}{2}y\der z+\der q)^2.
$$
This shows that it belongs to the Gibbons-Hawking class \cite{bi:gh} and 
that its self-dual part of the Weyl tensor vanishes.\\

\noindent
We also recall \cite{bi:przan2} that a suitable Lie-Backlund 
transformation brings equation (\ref{eq:k}) 
to the Boyer-Finley-Pleba\'nski \cite{bi:bf,bi:fp} equation
\footnote{also known to describe the {\bf SU}$(\infty)$ Toda lattice}
$$
F_{yy}+F_{zz}+({\rm e}^F)_{xx}=0
$$
for one real function $F=F(x,y,z)$ of three real variables. It is 
interesting to note that the metric (\ref{eq:metr}) corresponds to the 
simplest solution $F=0$ of this equation.\\

\noindent
Summing up we have the following theorem. 
\bt
Let $(z_1,\overline{z_1}, z_2,\overline{z_2})$ be coordinates on ${\cal U}\subset{\bf R}^4\cong {\bf C}^2$. The Riemannian manifold 
(${\cal U}'$, $g$), where 
$$
{\cal U}'
=\{{\cal U}\ni (z_1,z_2)\quad{\rm s.t.}\quad 
v-2z_2\overline{z_2} >0,\quad v=z_1+\overline{z_1}\}
$$
and 
$$
g=\frac{1}{(v-2z_2\overline{z_2})^{1/2}}
(\der z_1-2\overline{z_2}\der z_2)(\der\overline{z_1}-
2z_2\der\overline{z_2})+
4(v-2z_2\overline{z_2})^{1/2}\der z_2\der\overline{z_2},
$$
is Ricci-flat, anti-self-dual and has the anti-self-dual part of the Weyl 
tensor of type D. 
Moreover, (${\cal U}'$, $g$) admits a circle of almost-K\"ahler 
non-K\"ahler structures 
$$
J^+_{{\rm e}^{i\phi}}=2{\rm Re}\{
i{\rm e}^{i\phi}[\frac{1}{2(v-2z_2\overline{z_2})^{1/2}}
(\der z_1-2\overline{z_2}\der z_2)
\otimes 
(\partial_{\bar{2}}+2z_2\partial_{\bar{1}})-2
(v-2z_2\overline{z_2})^{1/2}\der z_2\otimes
\partial_{\bar{1}}]\}.
$$
These structures are 
parametrized by the real constant $\phi\in[0,2\pi[$. Their 
fundamental 2-forms are given by 
$$
\omega^+_{{\rm e}^{i\phi}}=i({\rm e}^{i\phi}\der z_2\dz\der z_1-
{\rm e}^{-i\phi}\der\overline{z_2}\dz\der\overline{z_1}). 
$$
\et

\noindent
5. Interestingly, our examples can be globalized.\\
Indeed, the transformation
$$
t=\frac{1}{2}\log (v-2z_2\overline{z_2}),\quad\quad\quad y=z_2+\overline{z_2},
\quad\quad\quad z=i(\overline{z_2}-z_2),\quad\quad\quad 
q=\frac{z_1-\overline{z_1}}{2i}
$$
brings the structures 
$(g, J^+_{{\rm e}^{i\phi}}, \omega^+_{{\rm e}^{i\phi}})$ of 
Theorem 1 to a form which is regular for all the values of the real 
parameters $(t,y,z,q)\in{\bf R}^4$.\\

\noindent
6. Finally, we observe that the metric (\ref{eq:metr}), as beeing 
anti-self-dual, possesses a strictly K\"ahler structure. This is given by 
$$
J=i[(\der z_1-2\overline{z_2}\der z_2)\otimes\partial_1-
(\der\overline{z_1}-2z_2\der\overline{z_2})\otimes\partial_{\bar{1}}+
\der\overline{z_2}\otimes (\partial_{\bar{2}}+2z_2\partial_{\bar{1}})-
\der z_2\otimes(\partial_2+2\overline{z_2}\partial_1)]
$$
and belongs to the structures of opposite orientation that 
$J^+_{{\rm e}^{i\phi}}$. 
It is interesting whether there exist Ricci-flat metrics that 
admit almost-K\"ahler non-K\"ahler structures but do not admit 
any strictly K\"ahler structure.\\

\noindent
{\bf Acknowledgements}\\
We are very grateful to W\l odek Jelonek for bringing to our attention the 
problem of existence of almost-K\"ahler non-K\"ahler Einstein metrics. We 
also wish to thank John Armstrong and Simon Salamon for information about 
their proof of existence of such metrics in four dimensions.\\
This work was completed during the workshop ``Spaces of geodesics and complex 
methods in general relativity and differential geometry'' held in Vienna at 
the Erwin Schr\"odinger Institute. We thank the members of the Institute and, 
especially, Helmuth Urbantke for creating the warm atmosphere during our 
stay in Vienna.

\end{document}